# AGOH'S CONJECTURE: ITS GENERALIZATIONS, ITS ANALOGUES

Andrei Vieru

**Abstract**

In this paper we formulate two generalizations of Agoh's conjecture. We also formulate conjectures involving congruence modulo primes about hyperbolic secant, hyperbolic tangent, Nörlund numbers, as well as about coefficients of expansions in powers of other analytic functions. We formulate a thesis about combinatorial objects that do not produce fake primes.

## 0. INTRODUCTION

One of the reasons Agoh's conjecture is said to be important lies in the fact that, if it is true – and it probably is – it permits to formulate a characterization of primes in terms of Bernoulli Numbers. It happens that Bernoulli Numbers are the best known, the most important and popular of a whole series of sequences of *combinatorial objects* – as they indeed may be thought of – that probably permit as well a characterization of primes in terms of these other combinatorial objects. Along with generalizations of Agoh's Conjecture in terms of Poly-Bernoulli Numbers, we'll try to provide a series of 'Agoh-like' conjectures involving other combinatorial objects than those aforementioned, which would be also liable – if they are true – to provide such characterizations of primes.

Our thesis is that the common features of these numbers – or objects – lie in the fact that they are all somehow related to $e$ and $\pi$.

It has to be noted that some other types of combinatorics – including counting of graph subsets with given properties (e.g. Perrin Sequence, Lucas sequence, generalizations of such sequences, etc.) always yield pseudo-primes. By the way they are not related neither to $e$ nor $\pi$, but they are indeed connected to algebraic integers, for example to some roots of characteristic polynomials in the case of linear recurrence sequences such as Perrin or Lucas or Fibonacci sequences.

# 1. BERNOULLI NUMBERS

The Bernoulli numbers $B_n$ are defined by the power series

$$\frac{ze^z}{e^z-1} = \sum_{n=0}^{\infty} B_n \frac{z^n}{n!}, \quad |z| < 2\pi$$

where all numbers $B_n$ are zero with odd index $n > 1$.

**CONJECTURE 1 (Generalized Agoh's conjecture)**

For any $q \geq 1$, for any odd $p \geq 3$

$pB_{q(p-1)} \equiv -1 \pmod{p} \Leftrightarrow p$ is prime

**GENERALIZED GIUGA'S CONJECTURE (JOSÉ MARIA GRAU, ANTONIO M. OLLER-MARCÉN)**

For any $q \geq 1$, for any odd $p \geq 1$

$1^{q(p-1)} + 2^{q(p-1)} + \ldots + (p-1)^{q(p-1)} \equiv -1 \pmod{p} \Leftrightarrow p$ is prime

In the particular case when $q = 1$, these conjectures are known, respectively, as Agoh's conjecture and Giuga's conjecture.

In the original formulation of the Agoh's conjecture there is no need to specify that $p$ is supposed to be odd, because $p - 1$ is supposed to be even.

In the formulation of the *Generalized Agoh's Conjecture* there is such a need, because otherwise $p$ and $q(p-1)$ might be both even and then the statement does not hold anymore. For $q = 4$ and $p = 30$, we find that

$30\, B_{4 \times 29} = 30 B_{116} \equiv -1 \pmod{30}$ **(I)**

(and, as it will become obvious later in this article, this congruence has something to do with the equalities $30 = 2 \times 3 \times 5$ and $2 \times 3 \times 5 - (2 \times 3 + 2 \times 5 + 3 \times 5) = -1$)

When formulating the *generalized Giuga's conjecture*, there also is such a need – i.e. to make it clear that $p$ is supposed to be odd – otherwise we find the following counterexample, with $q = 4$ and $p = 30$

$$\sum_{j=1}^{j=29} j^{4 \times 29} = \sum_{j=1}^{j=29} j^{116} =$$

44220744024334411808277004766010821764066504495575002284763550490669204193416663839965975095737376808664013068829912712275553732532548561543156543106336971952922034919 ≡ − 1 (mod 30)

In *The Equivalence of Giuga's and Agoh's Conjectures*, Bernd C. Kellner provided a proof of their equivalence.

Before we shall write a proof of the generalized Agoh's Conjecture, let us give an example:

$B_{60} = -1215233140483755572040304994079820246041491 / 56786730$

It is easy to check by direct computation that:

$61 B_{60} \equiv -1 \pmod{61}$ $\qquad$ $7 B_{60} \equiv -1 \pmod{7}$

$31 B_{60} \equiv -1 \pmod{31}$ $\qquad$ $5 B_{60} \equiv -1 \pmod{5}$

$13 B_{60} \equiv -1 \pmod{13}$ $\qquad$ $3 B_{60} \equiv -1 \pmod{3}$

$11 B_{60} \equiv -1 \pmod{11}$

## 2. POLY-BERNOULLI NUMBERS

Poly-Bernoulli numbers of order $k$, $B^{(k)}_n$, are defined by Masanobu Kaneko using polylogarithms of order $k$ in the following beautiful way. According to Kaneko:

For every integer $k$, a sequence of rational numbers is defined $B^{(k)}_n$ ($n = 0,1,2,...$), which is refered to as poly-Bernoulli numbers, by

$$\frac{1}{z} \operatorname{Li}_k(z) \Big|_{z=1-e^{-x}} = \sum_{n=0}^{\infty} B^{(k)}_n \frac{x^n}{n!}$$

Here, for any integer $k > 0$, $\operatorname{Li}_k(z)$ denotes the formal power series

$$\sum_{m=1}^{\infty} \frac{z^m}{m^k}$$

When $k = 1$, $B^{(1)}_n$ is the usual Bernoulli number (with $B^{(1)}_1 = \frac{1}{2}$):

$$\frac{xe^x}{e^x - 1} = \sum_{n=0}^{\infty} B^{(1)}_n \frac{x^n}{n!}$$

**CONJECTURE 2**

For any $k \geq 1$, for any $q \geq 1$ and for any odd $p > 1$

$p^k B^{(k)}_{q(p-1)} \equiv -1 \pmod{p} \Leftrightarrow p$ is prime

## 3. CAUCHY NUMBERS OF THE FIRST KIND

The Cauchy numbers $a_0, a_1, a_2, \ldots$ of the first kind are defined by

$$\sum_{n=0}^{\infty} a_n \frac{t^n}{n!} = \frac{t}{\ln(1+t)} \qquad |t| < 1$$

The first few are $a_0 = 1$, $a_1 = 1/2$, $a_2 = -1/6$, $a_3 = 1/4$, $a_4 = -19/30$, $a_5 = 9/4$, $a_6 = -863/84$, $a_7 = 1375/24$, etc.

They obey the following nice formula:

$$a_n = \sum_{j=0}^{n} \frac{s(n,j)}{j+1}$$

where $s(n, j)$ are the signed Stirling numbers of the first kind.

Let $a_n$ designate the $n$-th Bernoulli number of the second kind. Let $D(a_n)$ and $N(a_n)$ designate, respectively, its numerator and its denominator. Let $p(k)$ designate the $k$-th odd prime.

### CONJECTURE 3

$\forall P \geq 3, \forall m \geq 0, D(a_{2^m P}) \equiv N(a_{2^m P}) \pmod{P} \Leftrightarrow [P \text{ is prime} \wedge \exists k > m \text{ such as } P = p(k)]$

### CONJECTURE 4

For any odd primes $p_j$ and $p_k$, $D(a_{p_j p_k}) \equiv N(a_{p_j p_k}) \pmod{p_j} \Leftrightarrow p_j \geq p_k$

### CONJECTURE 5

$\dfrac{D(a_{2n})}{D(a_{2n+1})} = \lambda(n + \dfrac{1}{2})$ where $\lambda$ is an integer (most often 1) or a rational which, in decimal expansion, has at most one digit before the period.

### CONJECTURE 6

For any $q$ and for any odd $p > 2$, $p a_{q(p-1)} \equiv (-1)^{q-1} \pmod{p} \Leftrightarrow p$ is prime and $q \leq p$

### CONJECTURE 7

If $N_{2n+1}$ ($n > 3$) is the numerator of a Cauchy number of the first kind of the (not necessarily prime) odd rank $2n+1$, then it is divided by $(2n-1)^3$.

$\forall n > 5 \ [(2n-1)^3 \text{ divides } N_{2n+1} \wedge (2n-1)^4 \textit{ do not } \text{divide } N_{2n+1}] \Leftrightarrow 2n-1$ is prime

## 4. BERNOULLI NUMBERS OF THE SECOND KIND

The Bernoulli numbers of the second kind are defined in a slightly modified way:

$$\frac{t}{\ln(1+t)} = \sum_{n=0}^{\infty} b_n t^n \qquad |t|<1$$

Clearly, $b_n = \frac{a_n}{n!}$

So, the first few are :
$b_0 = 1$, $b_1 = 1/2$, $b_2 = -1/12$, $b_3 = 1/24$, $b_4 = -19/720$, $b_5 = 3/160$, $b_6 = -863/60480$, $b_7 = 275/24192$, etc.

$b_n$ obey the following beautiful formula:

$$\sum_{k=0}^{n} \frac{(-1)^k b_{n-k}}{k+1} = \delta_{n,0},$$

where $\delta_{n,0} = 1$ or 0 according to whether $n = 0$ or not.

### CONJECTURE 8
For any $p > 2$ $pb_p \equiv -1 \pmod{p} \Leftrightarrow p$ is prime

## 5. SIGNED TANGENT AND SECANT NUMBERS

Let us consider the coefficients of Taylor series of the sum of the hyperbolic tangent and hyperbolic secant, i. e. the numbers $S_n$ defined by

$$\tanh(z) + \text{sech}(z) = \sum_{n=0}^{\infty} S_n z^n$$

This yields, $S_0 = 1$, $S_1 = 1$, $S_2 = -1/2$, $S_3 = -1/3$, $S_4 = 5/24$, $S_5 = 2/15$, $S_6 = -61/720$, etc.

### CONJECTURE 9
(1) For any $p > 1$, $pS_p \equiv -1 \pmod{p} \Leftrightarrow p$ is prime
(2) For any $p > 2$, $pS_{p+1} \equiv 1 \pmod{p} \Leftrightarrow p$ is prime
(3) For any $p > 2$, $pS_{2p-1} \equiv 1 \pmod{p} \Leftrightarrow p$ is prime

Although it is well-known that the tangent and secant numbers are related to Bernoulli numbers (of the first kind) and to Euler numbers, it is rather surprising that they obey, in the same time, to the three statements of the conjecture. Some of the S numbers manage to show, at once, more than one congruence. For example, according to **(1)**, $5S_5 \equiv -1 \pmod{5}$ and, according to **(3)**, $3S_5 \equiv 1 \pmod{3}$.

According to **(1)**, $13S_{13} \equiv -1 \pmod{13}$ and, according to **(3)**, $7S_{13} \equiv 1 \pmod{7}$, etc., etc.

## 6. NÖRLUND NUMBERS (CAUCHY NUMBERS OF THE SECOND KIND)

The Nörlund numbers $a_n$ are defined by

$$\frac{t}{(1+t)\ln(1+t)} = \sum_{n=0}^{\infty} a_n \frac{t^n}{n!}$$

The first few, for $n = 0, 1, 2\ldots$, are:

$1, -1/2, 5/6, -9/4, 251/30, -475/12, 19087/84, -36799/24, 1070017/90,$

**CONJECTURE 10**

For any $q$ and for any odd $p > 2$, $pN_{q(p-1)} \equiv (-1)^{q-1} \pmod{p} \Leftrightarrow p$ is prime and $q \leq p$

**CONJECTURE 11**

The numerator of any Nörlund number of odd rank (not necessarily prime) is divided by the square power of its rank.

It is also divided by the cube power (or higher) of its own rank if and only if this rank *is not* prime.

If the odd rank of a Nörlund number has only two prime factors, $p_1$ and $p_2$ ($p_1 < p_2$), then $(p_1 p_2)^{p_1+1}$ divides its numerator.

## 8. ERRATUM

The proof of Agoh-Giuga Conjecture seems to be a desperate case.
In a previous version of this paper we made an erroneous attempt of a proof, that contained an mistake (pointed out by Professor Robin Chapman to whom we express our deepest gratitude).

## 9. ON A PROPERTY OF THE THIRTEEN KNOWN GIUGA NUMBERS

I wonder if the following remark might be useful:

Among the two biggest prime divisors of all thirteen already known Giuga numbers there is exactly one prime of the form $4k+1$, while all other odd prime divisors of the aforementioned numbers are of the form $4k+3$. Moreover, if P is an odd prime divisor of some already discovered Giuga number G, then either

1) P is of the form $4k+3$ and then $(G/P-1)/P$ is also of the form $4k+3$ (although not necessarily prime)

or

2) P is of the form $4k+1$ and then $(G/P-1)/P$ is also of the form $4k+1$ (although not necessarily prime)

Should we conjecture that all (even) Giuga numbers display the same feature?

Lava's Conjecture has recently given raise to scepticism and criticism, because thirteen examples are usually considered as being *not enough*.

Probably, searching new Giuga numbers in this direction would soon lead to discoveries of yet unknown numbers of this kind. However, they wouldn't underpin any conjecture because the very method of search would contain a *petitio principii.*

## 10. OTHER AGOH-LIKE CONJECTURES

Many integer sequences that are coefficients of analytic functions expansions (conjecturally) show congruence (mod $p$) to ±1 for terms of rank $p$ or $p - 1 \Leftrightarrow p$ is prime. Examples of conjectures involving congruence modulo primes:

In the sequence OEIS A002105 (Reduced tangent numbers)

$\forall p > 0 \quad A_p \equiv 1 \pmod{p} \Leftrightarrow p$ is prime  (3°)

In A064062 Generalized Catalan numbers

$A_p \equiv 1 \pmod{p} \Leftrightarrow p$ is prime  (4°)

In A000587   Rao Uppuluri-Carpenter numbers

$A_{p+2} \equiv 1 \pmod{p} \Leftrightarrow p$ is prime  **(5°)**

In A001586   Generalized Euler Numbers (or Springer Numbers)

$A_p \equiv 1 \pmod{p} \Leftrightarrow p$ is a prime of the form $4n+1$

or $p$ is a power of 2  **(6°a)**

$A_p \equiv -1 \pmod{p} \Leftrightarrow p$ is a prime of the form $4n-1$  **(6°b)**

for any $n \geq 1$  $A_{3^n} \equiv 2 \cdot (-1)^{n+1} \pmod{3^n}$

In A000111 Euler up/down numbers

$A_p \equiv 1 \pmod{p} \Leftrightarrow p$ is a prime of the form $4n+1$

or $p$ is a power of 2  **(7°a)**

$A_p \equiv -1 \pmod{p} \Leftrightarrow p$ is prime of the form $4n-1$  **(7°b)**

for any $n \geq 1$  $A_{3^n} \equiv 2 \cdot (-1)^{n+1} \pmod{3^n}$

In A007836  Springer numbers associated with symplectic group

$A_p \equiv 1 \pmod{p} \Leftrightarrow p$ is a prime of the form $4n+1$  **(8°a)**

$A_p \equiv -1 \pmod{p} \Leftrightarrow p$ is a prime of the form $4n-1$  **(8°b)**

or $p$ is a power of 2

In A001006 Motzkin numbers

For any $p > 1$   $M_p \equiv 1 \pmod{p} \Leftrightarrow p$ is a prime or the square of a prime

In A000108 Catalan numbers

For an enough large value of $p$, $C(p-1) \equiv -1 \pmod{p} \Leftrightarrow p$ is a prime or the square of a prime.

**CONCLUSION**

Surely, there are lots of such sequences of numbers – we prefer saying "combinatorial objects" – for whom Agoh-like Conjectures might be formulated. The question is whether the thesis according to which combinatorial objects related to $e$ or/and $\pi$ do not produce fake primes holds or not. In particular we strongly "believe" in the converse of the Wolstenholme's Theorem, as far as Harmonic Numbers satisfy the criteria our thesis.

Andrei Vieru
andreivieru007@gmail.com
Homepage: www.andreivieru.com